\documentclass[english,12pt]{amsart}
\usepackage{amsmath}
\usepackage{amsthm}
\usepackage{pdfsync}
\usepackage{a4wide}
\usepackage{amssymb}
\usepackage{hyperref}
\usepackage{amsfonts}
\usepackage{color}
\usepackage{enumerate}
\usepackage{graphicx}
\usepackage{amsfonts}
\usepackage{color}
\usepackage{enumerate}
\theoremstyle{definition}

\def\eps{\epsilon}
\def\rr{\mathbb{R}}

\def\intr{\int_{\rr^N}}
\def\rr{\mathbb{R}}
\def\intr{\int_{\rr^N}}

\def\n{\nabla}
\def\nn{\mathbb{N}}
\def\eps{\epsilon}

\numberwithin{cor}{section} \numberwithin{lema}{section}
\numberwithin{prop}{section} \numberwithin{Rem}{section}
\numberwithin{conj}{section} 
\numberwithin{Def}{section}
\newtheorem{theorem}{Theorem}[section]

\theoremstyle{def}

\newtheorem{remark}[theorem]{Remark}

\usepackage[normalem]{ulem}
\definecolor{DarkBlue}{rgb}{0,0.1,0.7} 
\definecolor{DarkGreen}{rgb}{0,0.5,0.1} 

\newcommand\soutD{\bgroup\markoverwith
	{\textcolor{DarkGreen}{\rule[.5ex]{2pt}{1pt}}}\ULon}
\newcommand\soutP{\bgroup\markoverwith
	{\textcolor{blue}{\rule[.5ex]{2pt}{1pt}}}\ULon}
\newcommand{\Hm}[1]{\leavevmode{\marginpar{\tiny%
			$\hbox to 0mm{\hspace*{-0.5mm}$\leftarrow$\hss}%
			\vcenter{\vrule depth 0.1mm height 0.1mm width \the\marginparwidth}%
			\hbox to
			0mm{\hss$\rightarrow$\hspace*{-0.5mm}}$\\\relax\raggedright #1}}}

\begin{document}
	
	\title[Sharp weighted Hardy-Rellich inequalities]{Weighted Hardy-Rellich type inequalities: improved best constants and  symmetry breaking}
\author{Cristian Cazacu}
\address{Cristian Cazacu: Faculty of Mathematics and Computer Science \\
	University of Bucharest\\
	14 Academiei Street \\
	010014 Bucharest, Romania
	\&
	Gheorghe Mihoc-Caius Iacob Institute of Mathematical
	Statistics and Applied Mathematics of the Romanian Academy\\ No.13 Calea 13 Septembrie, Sector 5\\
	050711 Bucharest, Romania}
\email{cristian.cazacu@fmi.unibuc.ro}
		
		\medskip 
\author{Irina Fidel}
\address{Irina Fidel: Sfantul Vasile Secondary School \\ No. 145 Republicii\\
		100389 Ploiesti, Romania\\ 
			\&	
		Faculty of Mathematics and Computer Science \\
		University of Bucharest\\
		14 Academiei Street \\
		010014 Bucharest, Romania}
		\email{irina.fidel@s.unibuc.ro}
	

	\begin{abstract}
		When studying the weighted Hardy-Rellich inequality in $L^2$ with the full gradient replaced by the radial derivative    the best constant becomes trivially larger or equal than in the first situation. Our contribution is to determine the new sharp constant and to show that for some part of the weights is strictly larger than before. In some cases we emphasize that the extremals functions of the sharp constant are not radially symmetric.  
	
	\end{abstract}
\maketitle

\section{Introduction}	
 
The celebrated $L^2$ \textit{Hardy inequality} (e.g. \cite{HLP, OK}) states that for $N\geq 3$ and $u\in C_c^\infty(\rr^N)$ it holds  
\begin{equation}\label{HI}
	\intr|\nabla u|^2dx\geq C_H\intr \frac{|u|^2}{|x|^2}dx, \quad C_H(N):=\frac{(N-2)^2}{4},
\end{equation}
where the constant $C_H(N)$ is sharp.

\textit{Rellich inequality} (e.g. \cite{R}) asserts  that for $N\geq 5$ and $u\in C_c^{\infty}(\mathbb{R}^N)$ we have  
\begin{equation}\label{RI}
\intr|\Delta u|^2dx\geq C_R(N)\intr \frac{|u|^2}{|x|^4}dx, \quad C_R(N):=\left(\frac{N(N-4)}{4}\right)^2,
\end{equation} with the best constant $C_R(N)$.  The Hardy and Rellich inequalities are important tools widely used in the analysis of partial differential operators and equations of harmonic and biharmonic-type.

The \textit{Hardy-Rellich inequality} has been studied more recently (see, e.g.  \cite{Tertikas, Ghoussoub, Cazacu}). This is in fact an improved Hardy inequality (with a larger optimal constant) applied to classes of vector fields originated from potential gradients, which arises in fluid mechanics.  For $N\geq 3$ and $u\in C_0^{\infty}(\mathbb{R}^N)$ this leads to  
\begin{equation}\label{HRI}
\intr|\Delta u|^2dx\geq C_{HR}(N)\intr\frac{|\nabla u|^2}{|x|^2}dx,
\end{equation}
 with the best constant 
$$C_{HR}(N):=\left\{
\begin{aligned}
	&\frac{N^2}{4}, & N\geq 5\\
	&3, &N=4\\
	&\frac{25}{36},& N=3.
\end{aligned}
\right.  $$

 The Hardy-Rellich inequality was firstly extended in \cite{Tertikas} to more general singular weights of the form 
\begin{equation}\label{weight_ineq1}
	\intr |\Delta u|^2|x|^{m} dx \geq C(N, m) \intr |\nabla  u|^2|x|^{m-2} dx, \quad \forall u\in C_c^\infty(\rr^N),
\end{equation} 
where the authors proved that for any $N\geq 5$ and any $4-N < m\leq 0$ the best constant is 
\begin{equation}\label{best_const_TZ}
	C(N, m)=\min_{k=0,1,2...}\frac{\left(\frac{(N-4+m)(N-m)}{4}+k(N+k-2)\right)^2}{\left(\frac{N-4+m}{2}\right)^2+k(N+k-2)}.
\end{equation}
In particular, according to the computations in \cite[pag. 453]{Tertikas}, if $\frac{N+4-2\sqrt{N^2-N+1}}{3}\leq m\leq 0$ then 
$$C(N, m)=\left(\frac{N-m}{2}\right)^2$$ 
otherwise, if $4-N< m< \frac{N+4-2\sqrt{N^2-N+1}}{3}$ then 
$$C(N, m)< \left(\frac{N-m}{2}\right)^2.$$

To ensure the integrability of the singular term in inequality \eqref{weight_ineq1} we need to impose that $|x|^{m-2}\in L_{loc}^1(\rr^N)$ which gives us the constraint 
\begin{equation}\label{constraint}
	m> 2-N \quad (or \ \ m+N-2>0).
\end{equation}

The weighted inequality  \eqref{weight_ineq1}  was later extended in \cite{Ghoussoub} to all the cases $N\geq 1$ and $m >2-N$. Optimal constants of the cases which were not covered in \cite{Tertikas} were solved in \cite[Theorem 6.1]{Ghoussoub}. Next, we emphasize a brief presentation of these additional cases: 
\begin{itemize}
\item If $N=1$ and $ m\in\left(1,\frac{7}{3}\right] \cup \left[3, \infty\right)$ then $C(1, m)=\left(\frac{1-m}{2}\right)^2$.

\item If $N=1$ and $m\in(\frac{7}{3}; 3)$ then $C(1,m)\leq \left(\frac{1-m}{2}\right)^2$ (there are values $m$ for which the inequality is strict).
\item If $N\geq 1$ and $m=4-N$ then $C(N,m)= \min \{(N-2)^2, N-1\}$.
\item If $N\geq 2$ and $ \frac{N+4-2\sqrt{N^2-N+1}}{3}\leq m$ then  $C(N,m)=\left(\frac{N-m}{2}\right)^2$.
\item If $2\leq N\leq 3$ and $2-N< m < \frac{N+4-2\sqrt{N^2-N+1}}{3}$, or $N\geq 4$ and $2-N< m\leq 4-N$, then $C(N,m)=\frac{\left(\frac{(N-4+m)(N+m)}{4}+N-1\right)^2}{\left(\frac{N-4+m}{2}\right)^2+N-1}$
\item If $N=3$ and $m \geq  \frac{N+4-2\sqrt{N^2-N+1}}{3}$ or $N\geq 4$ and $m> 4-N$ the best constant requires a further subdivision,  based on very technical  expressions, which could be consulted in \cite[Th. 6.1, pag. 52]{Ghoussoub}. 
\end{itemize}
Subsequent extensions of the weighted Hardy-Rellich type inequalities with reminder terms have been done recently in \cite{Tertikas2} and \cite{J} by applying factorization methods. Also recent improvements when adding magnetic fields have been established in \cite{CCF, LG}. See also very  recent results on the Hardy-Rellich inequalities in \cite{B, MR} and in \cite{Ham} (in the context of solenoidal vector fields) and the references there in.   

Overall for any $N\geq 1$ and $m>2-N$ always happens that the best constant in \eqref{weight_ineq1} does not pass the threshold
\begin{equation}\label{threshold}
	C(N, m)\leq \left(\frac{N-m}{2}\right)^2.
\end{equation} 

In this paper we study a new weighted Hardy-Rellich type inequality which, to our knowledge, has not been treated yet in the literature, by replacing the full gradient in \eqref{weight_ineq1} with the radial derivative $\partial_r u:=\frac{x}{|x|}\cdot \nabla u$ on the right hand side, namely 
\begin{equation}\label{newHR}
	\intr |\Delta u|^2|x|^{m} dx \geq \tilde{C}(N, m) \intr |x \cdot \nabla  u|^2|x|^{m-4} dx, \quad \forall u\in C_c^\infty(\rr^N). 
\end{equation}
where $\tilde{C}(N, m)$ denotes the best constant in \eqref{newHR}.  Since $|\partial_r u|\leq |\nabla u|$ notice that inequality \eqref{newHR} is a relaxation of \eqref{weight_ineq1} and in balance with that the best constant might be larger or equal, i.e. $C(N, m)\leq \tilde{C}(N, m)$. 

Our purposes are to supply explicitly the best constant $\tilde{C}(N, m)$ for the full range of parameters $N\geq 1$ and $m > 2-N$, to emphasize situations in which we get an improvement in \eqref{newHR} with respect to \eqref{weight_ineq1}  i.e. 
 $\tilde{C}(N, m)> C(N, m)$ and to put in evidence the radially symmetry breaking for the optimal approximations of the sharp constants. 

\section{Main result}

On the purpose to state the main result we need to introduce a cut-off  function $g\in C^{\infty}([0,\infty))$, with $0\leq g\leq 1$,  such that 
\begin{equation}\label{function_g} 
	g(r)=\left\{
	\begin{array}{ll}
		1, & 0\leq r\leq 1,\\
		0, & r\geq 2.\\
	\end{array}
	\right. 
\end{equation}

The main result of this paper is the following 
\begin{theorem}\label{key1} Let $N\geq 1$ and $m >2-N$. Then  
	\begin{equation}\label{weight_ineq29}
		\intr |\Delta u|^2|x|^{m} dx \geq \tilde{C}(N, m) \intr |x \cdot \nabla  u|^2|x|^{m-4} dx. \quad \forall u\in C_c^\infty(\rr^N), 
	\end{equation}
	where the optimal constant $\tilde{C}(N, m)$ is given as folllows.
	


	If $N=1$ and $m>1$ or if $N\geq 2$ and $m\in  [2-\sqrt{(N-1)^2+1},  2+\sqrt{(N-1)^2+1}]$, then $$\tilde{C}(N,m)=\left(\frac{N-m}{2}\right)^2,$$ which is approximated by the sequence 
$\{u_\eps\}_{\eps>0}$ given by 
\begin{equation}\label{min_seq_2}
	u_\eps (x) =|x|^{-\frac{N+m-4}{2}+\eps}g(|x|).
\end{equation}

		If $N\geq 2$ and $m\in (2-N, 2-\sqrt{(N-1)^2+1})$, then $$\tilde{C}(N,m)= \frac{\left((m-2)^2-N^2\right)^2} {4(N+m-4)^2}$$ which is approximated by the  sequence 
	$\{u_\eps\}_{\eps>0}$ given by 
	\begin{equation}\label{min_seq_22}
		u_\eps (x) =|x|^{-\frac{N+m-4}{2}+\eps}g(|x|)\phi_1(x),
	\end{equation}
where $\phi_1$ is a spherical harmonic function of degree $1$ with  $ \|{\phi_1}\|_{L^2(S^{N-1})}=1$.

		If $N\geq 2$ and $m\in (2+\sqrt{(N-1)^2+1}, \infty)$, then $$\tilde{C}(N,m)=\min_{l\leq k(m)} \tilde{C}(N, m, l),\quad \tilde{C}(N, m, l):=\frac{\left(-N+m-2l\right)^2\left(2l+m+N-4 \right)^2}{4\left(m+N-4\right)^2}$$ where $k(m)$ is defined later in \eqref{ikmn}-\eqref{k's}. The constant $\tilde{C}(N, m)$ is approximated by the sequence 
		given by 
		\begin{equation}\label{min_seq_21}
			u_\eps (x) =|x|^{-\frac{N+m-4}{2}+\eps}g(|x|)\phi_{l_{\min}}(x),
		\end{equation}
		where $\phi_{l_{\min}}$ is a spherical harmonic function of degree $l_{\min}$ such that $\|\phi_{l_{\min}}\|_{L^2(S^{N-1})}=1$
	and $l_{\min}:=\arg\min_{l\leq k(m)}\tilde{C}(N, m, l)$. 
\end{theorem} 
\begin{remark}
	Notice also that in all the situations above $$\tilde{C}(N, m)\leq \left(\frac{N-m}{2}\right)^2,$$
	but there are cases (see for instance $N=1$ and $m\in (\frac{7}{3}, 3)$) when our best constant $\tilde{C}(N, m)$ in \eqref{newHR} improves with respect to the best constant $C(N, m)$ in \eqref{weight_ineq1}.  
\end{remark}
\begin{remark} The approximating sequences $u_\epsilon$ in \eqref{min_seq_2}, \eqref{min_seq_21}, \eqref{min_seq_22} do not belong to the space $C_c^\infty(\rr^N)$ but they are in the energy space of the inequality \eqref{weight_ineq29}, i.e. both terms in \eqref{weight_ineq29} are finite for $u_\epsilon$. So, by regularizing $u_\epsilon$ near the origin one can show that the constants remain sharp for functions $u\in C_c^\infty(\rr^N)$, see for instance similar arguments in \cite{CFL, CIM}.    
\end{remark}
\section{Proof of the main result} 
	We will use spherical coordinates instead of cartesian coordinates. This coordinates transformation is given by  
	$$ x\in \rr^N\setminus\{0\} \mapsto (r, \sigma) \in (0, \infty)\times S^{N-1}, \quad r=|x|, \quad \sigma=\frac{x}{|x|},$$ where $S^{N-1}$ is the $N-1$-dimensional sphere with respect to the Hausdorff measure in $\rr^N$. We will use the formula of the Laplacian in spherical coordinates
	\begin{equation}\label{Laplacian}
		\Delta=\partial^2_{rr}+\frac{N-1}{r} \partial_r+\frac{1}{r^2} \Delta_{S^{N-1}},
	\end{equation}
	where $\partial_r$ and $\partial^2_{rr}$ are first and second order partial derivatives  with respect to the radial component $r$ whereas (for fixed $r$) the Laplace-Beltrami operator with respect to the metric tensor  on  $S^{N-1}$ and respectively the spherical gradient are given by 
	$$\Delta_{S^{N-1}}u(r\sigma)=\Delta\left[u\left(\frac{x}{|x|}\right)\right]_{|x=\sigma}, \quad \nabla_{S^{N-1}}u(r\sigma)=\nabla\left[u\left(\frac{x}{|x|}\right)\right]_{|x=\sigma}.$$
	Applying the spherical harmonics decomposition we  can expand  $u\in C_c^\infty(\rr^N)$ as
	$$u(x)=u(r\sigma)=\sum_{k=0}^{\infty} u_k(r)\phi_k(\sigma).$$
 The set of functions $\{\phi_k\}_{k\geq 0}$ are spherical harmonics of degree $k$ which consists in an orthogonal basis in $L^2(S^{N-1})$.
	These functions satisfy the properties
	\begin{equation}\label{orthoeigen}
		\left\{
		\begin{array}{ll}
			-\Delta_{S^{N-1}}\phi_k=c_k \phi_k  \textrm{ on } S^{N-1}, \\[6pt]
			-\int_{S^{N-1}} \Delta_{S^{N-1}}\phi_k \phi_l d\sigma= \int_{S^{N-1}} \n_{S^{N-1}}\phi_k\cdot \n_{S^{N-1}} \phi_l d\sigma\\[3pt]
			=c_k \int_{S^{N-1}} \phi_k \phi_l d\sigma=c_k \delta_{lk},  \quad k, l\in \nn,
		\end{array}
		\right.
	\end{equation}
	where $c_k=k(k+N-2)$, $k\geq 0$ are the eigenvalues of the Laplace-Beltrami operator $\Delta_{S^{N-1}}$,
	where $\delta_{lk}$ represents the Kronecker symbol, see e.g. \cite{Groemer} for more detailed properties of spherical harmonics.

%
	Since in view of \eqref{Laplacian} $$\Delta\left[u_k(|x|)\phi_k\left(\frac{x}{|x|}\right)\right]=\left(\Delta_r u_k(r)-\frac{c_k}{r^2}u_k(r)\right)\phi_k(\sigma),$$
due to \eqref{orthoeigen} similar computations as in \cite{Cazacu} lead to 	
	\begin{equation}\label{Bil}
		\intr |\Delta u|^2|x|^m dx=\sum_{k=0}^{\infty} \int_0^\infty \left(|\Delta_r u_k|^2 +\frac{c_k^2}{r^4}u_k^2 -\frac{2c_k}{r^2} u_k \Delta_r u_k  \right) r^{N+m-1} dr.
	\end{equation}
		Next, we will write $u_k^{\prime}$ and $u_{k}^{\prime\prime}$ to express both first and second derivatives of the Fourier coefficients $\{u_k\}_k$. 
	
		We need to compute $\int_0^{\infty}|\Delta_r u_k|r^{N+m-1}$ and $ \int_0^{\infty}\frac{\Delta_r u_k}{r^2} u_k(r)r^{N+m-1}dr$. 
	
	Integration by parts leads to 
	\begin{equation}\label{deltauk}
		\begin{aligned}
			\int_0^\infty |\Delta_r u_k|^2 r^{N+m-1} dr &=(N-1)(1-m)\int_0^\infty |u_k^{\prime}(r)|r^{N+m-3} dr +\\
			&+\int_0^\infty |u_k^{\prime\prime }(r)|r^{N+m-1} dr,
		\end{aligned}
	\end{equation}
	and 
	\begin{equation}\label{deltauk2}
		\begin{aligned}
			\int_0^\infty \frac{\Delta_r u_k}{r^2} u_k(r) r^{N+m-1} dr&=\frac{1}{2}\left( -2 \int_0^\infty |u_k^{\prime}(r)|r^{N+m-3} dr\right) +\\
			&+(N+m-4)(m-2)\int_0^\infty |u_k(r)|r^{N+m-5} dr.
		\end{aligned}
	\end{equation}
	
Then 	\eqref{Bil} becomes 
	\begin{align}\label{form1}
		\int_{\rr^N} |\Delta u|^2|x|^m dx &=  \sum_{k=0}^{\infty} \Big(\int_0^\infty |u_k^{\prime\prime}|^2 r^{N+m-1} dr +\\
\nonumber &+ (2c_k+(N-1)(1-m))\int_0^\infty |u_k^{\prime}|^2 r^{N+m-3} dr\\
		\nonumber
		&+ (c_k^2-c_k(m-2)(N+m-4))\int_0^\infty |u_k|^2 r^{N+m-5} dr\Big). \nonumber
	\end{align}
	and more easily, since $\partial_ru=\frac{x}{|x|}\cdot \nabla u$, we get 
\begin{equation}\label{form21}
		\intr |x\cdot \n  u|^2|x|^{m-4} dx = \sum_{k=0}^{\infty} \int_0^\infty  |u_k^{\prime}|^2 r^{N+m-3} dr.
	\end{equation}

	\subsection{The case $N=1$ and $m>1$}

The inequality \eqref{weight_ineq29} reduces to:
	\begin{equation}\label{1d_ineq}
		\int_\rr |u^{\prime\prime}|^2 r^{m} dr \geq \tilde{C}(1, m)\int_\rr |u^\prime|^2 r^{m-2} dr, \quad \forall u \in C_c^\infty(\rr).   
	\end{equation} 
	
	The relation \eqref{1d_ineq} comes from 
	\begin{equation}
		\begin{aligned}
			\int_{\mathbb{R}}(u^\prime)^2 r^{m-2}dr&=\int_{\mathbb{R}}(u^\prime)^2\left(\frac{r^{m-1}}{m-1}\right)^\prime dr\\
&=-\frac{2}{m-1}\int_{\mathbb{R}}u^\prime u^{\prime\prime} r^{m-1}dr\\
			&=-\frac{2}{m-1}\int_{\mathbb{R}}u^\prime r^{\frac{m}{2}} u^{\prime\prime} r^{\frac{m-2}{2}}dr\\
			&\leq -\frac{2}{m-1}\left(\int_{\mathbb{R}}(u^{\prime\prime})^2 r^m dr \right)^{\frac{1}{2}}\left(\int_{\mathbb{R}}(u^\prime)^2 r^{m-2}dr \right)^{\frac{1}{2}}.
			\nonumber
		\end{aligned}
	\end{equation}
	
	Integration by parts and Cauchy-Schwarz inequality lead that \eqref{1d_ineq} hold with the constant 
	\begin{equation}\label{opt1}
	\tilde{C}(1, m)\geq \left(\frac{m-1}{2}\right)^2.
	\end{equation}
	\subsection{The case $N\geq 2$.} 

In this case it remains to compare the right hand sides in \eqref{form1}-\eqref{form21}. 
	We will apply the well-known 1-d weighted Hardy inequalities (see, e.g. \cite[page 454]{Tertikas}, \cite[Prop. 2.4]{Cassano})
	\begin{equation}\label{WR_1d_11}
		\int_0^\infty  |u_k^{\prime \prime}|^2  r^{N+m-1}dr \geq \left(\frac{N+m-2}{2}\right)^2  \int_0^\infty  |u_k^{\prime}|^2 r^{N+m-3} dr, 
	\end{equation}
	\begin{equation}\label{WR_1d_22}
		\int_0^\infty  |u_k^{\prime}|^2  r^{N+m-3}dr \geq \left(\frac{N+m-4}{2}\right)^2  \int_0^\infty  |u_k|^2 r^{N+m-5} dr. 
	\end{equation}
	\subsubsection{The case A: $(2-m)(N+m-4)\geq 0$}
	We distinguish two sub-cases as follows. 
	\medskip\\
		\textbf{A1)  $2-m\geq 0$ and $N+m-4\geq 0$.}      These imply  $4-N\leq m\leq 2$. 
		\medskip\\
 \textbf{A2) $2-m\leq 0$ and $N+m-4\leq 0$.}  Since $N\geq 2$ it follows that $m=N=2.$\\

Gluing both situations we can summarize that case {\bf A}   is equivalent to
\begin{equation}\label{caseA}
	 4-N\leq m\leq 2.
\end{equation} 

Then we obtain from \eqref{form1}-\eqref{form21} and \eqref{WR_1d_11} that 
	\begin{equation}\label{identi1}
		\begin{aligned}
			\intr |\Delta u|^2|x|^m dx&\geq \sum_{k=0}^{\infty}\left(\left(\frac{N+m-2}{2}\right)^2\int_{ 0}^{\infty}|u_k^{\prime}|^2 r^{N+m-3}dr\right.\\
			&\left.+(2c_k+(N-1)(1-m))\int_0^{\infty }|u_k^{\prime}|^2 r^{N+m-3}dr\right)\\
			&\geq \left(\left(\frac{N+m-2}{2}\right)^2+(N-1)(1-m)\right)\times\\
			&\times \sum_{k=0}^{\infty} \int_0^\infty  |u_k^{\prime}|^2 r^{N+m-3} dr\\
			&=\left(\frac{N-m}{2}\right)^2 \sum_{k=0}^{\infty}\int_0^{\infty}|u_k^{\prime}|^2 r^{N+m-3}dr\\
			&=\left(\frac{N-m}{2}\right)^2 \intr |x \cdot \nabla  u|^2|x|^{m-4} dx.\\
		\end{aligned}
	\end{equation}
	Hence 
	\begin{equation}\label{opt2}
		\tilde{C}(N, m)\geq \left(\frac{N-m}{2}\right)^2.
		\end{equation}
	
	\subsubsection*{The case B:  $(2-m)(N+m-4)< 0$.}	
\noindent	\textbf{B1) $2-m>0$ and $N+m-4<0$.} These combined with \eqref{constraint} are equivalent to 
	\begin{equation}\label{caseB1}
		2-N< m< 4-N.
	\end{equation}
	\textbf{B2) $m-2<0$ and $N+m-4>0$.} These are equivalent to
	 \begin{equation}\label{caseB2}
	 	m>2.
	 \end{equation} 
	
We first look at the spherical part (which contains the terms multiplied with $c_k$) in  the relation \eqref{form1}. Taking into account \eqref{WR_1d_22} and the fact that $c_k\geq 0$ for any $k\geq 0$  we successively have  
	\begin{align}\label{important1}
		2c_k\int_0^\infty &|u_k^{\prime}|^2 r^{N+m-3} dr + (c_k^2-c_k(m-2)(N+m-4))\int_0^\infty |u_k|^2 r^{N+m-5} dr\nonumber\\
		& \geq c_k \left(2\left(\frac{N+m-4}{2}\right)^2+c_k -(m-2)(N+m-4)\right) \int_0^\infty |u_k|^2 r^{N+m-5} dr
	\end{align}
	In view of this let us denote 
	\begin{equation}\label{ikmn}
I_k(m, N):=2\left(\frac{N+m-4}{2}\right)^2+c_k -(m-2)(N+m-4).\end{equation}

	In view of identity \eqref{form1}, Hardy inequalities \eqref{WR_1d_11}-\eqref{WR_1d_22} and \eqref{important1}-\eqref{ikmn} we get 
\begin{align}\label{formB2}
	\int_{\rr^N} |\Delta u|^2|x|^m dx &\geq  \left(\frac{N-m}{2}\right)^2\sum_{k=0}^{\infty}\int_0^\infty|u_k^{\prime}|^2 r^{N+m-3} dr+\nonumber\\
	&+ \sum_{k=0}^{\infty}c_k I_k(m,N)\int_0^\infty |u_k|^2 r^{N+m-5} dr.
\end{align}

	Next we want to investigate when $I_k(m,N)\geq 0$. We start with estimating 
	\begin{align}
I_1(m,N)&=
2\left(\frac{N+m-4}{2}\right)^2+N-1-(m-2)(N+m-4)\nonumber\\
&= -\frac 12 m^2 +2m +\frac 12 N^2-N-1.\nonumber
\end{align}

	The equation $I_1(m ,N)=0$ in the unknown $m$ has  the discriminant $\Delta=(N-1)^2+1\geq 0$ and  the roots  
	$$m_{1,2}=2\pm \sqrt{(N-1)^2+1}$$
	which imply 
	$$I_1(m, N)\geq 0, \quad \textrm{ iff } m \in [2-\sqrt{(N-1)^2+1},  \ 2+\sqrt{(N-1)^2+1}].$$
	Therefore, we conclude that 
	in the case {\bf B1)} we have 
	$$I_1(m, N)\geq 0, \quad \textrm{ iff } m \in [2-\sqrt{(N-1)^2+1},  \ 4-N)$$
	whereas in the case {\bf B2)}  we have 
	$$I_1(m, N)\geq 0, \quad \textrm{ iff } m \in (2, 2+\sqrt{(N-1)^2+1}].$$
	
	Since $\{I_k\}$ is an increasing sequence with respect to $k$ we get  that 
	\begin{equation}
		I_k(m,N)\geq 0, \quad \forall k\geq 1, 
	\end{equation} 
for any $m$ satisfying
\begin{equation}\label{conditions for m}
	 m\in [2-\sqrt{(N-1)^2+1},  \ 4-N) \cup (2, 2+\sqrt{(N-1)^2+1}].
\end{equation}
This together with \eqref{form1} and \eqref{formB2} yield to 
	\begin{equation}\label{opt3}
	\tilde{C}(N, m)\geq \left(\frac{N-m}{2}\right)^2, \quad \forall m \textrm{ as in } \eqref{conditions for m}. 
\end{equation}	

It remains to analyze the complementary "bad cases" of \eqref{conditions for m} for which $I_1(m, N)<0$:
	\begin{equation}\label{cazurile rele pt m}
		m\in \underbrace{(2-N, 2-\sqrt{(N-1)^2+1})}_{\textrm{remaining cases of B1)}}\cup \underbrace{(2+\sqrt{(N-1)^2+1}, \infty)}_{\textrm{remaining cases of B2)}}.
	\end{equation}
If $k=2$ we obtain 
\begin{align}
I_2(m, N)&=2\left(\frac{N+m-4}{2}\right)^2+2N -(m-2)(N+m-4)\nonumber\\
&=-\frac{m^2}{2}+2m+\frac{N^2}{2}.\nonumber
\end{align} Therefore
\begin{equation}\label{I_2 pozitiv}
I_2(m, N)\geq 0, \quad \textrm{ iff } m \in [2-\sqrt{N^2+4},  \ 2+\sqrt{N^2+4}].
	\end{equation}
This leads to 
\begin{equation}\label{I_k pozitiv}
	I_k(m,N)\geq 0, \quad \forall k\geq 2, \quad  m\in (2-N, 2-\sqrt{(N-1)^2+1}),
\end{equation} 
which cover the remaining cases of $\bold{B1)}.$


	\paragraph{Inequality \eqref{weight_ineq29} in the remaining cases of B1): $m\in (2-N, 2-\sqrt{(N-1)^2+1})$.}

 The right hand side in \eqref{form1} can be bounded from below in terms of a parameter $\varepsilon>0$ (which will be well precised later) as follows, in view of \eqref{formB2} and \eqref{WR_1d_22}:
\begin{align}\label{form100}
	\int_{\rr^N} |\Delta u|^2|x|^m dx &\geq \sum_{k=0, k\ne 1}^{\infty} \Bigg(\left(\frac{N-m}{2}\right)^2\int_0^\infty|u_k^{\prime}|^2 r^{N+m-3} dr+\nonumber\\
	&+ c_k I_k(m,N)\int_0^\infty |u_k|^2 r^{N+m-5} dr\Bigg)+\nonumber\\
	&+\left(\frac{N-m}{2}\right)^2\int_0^\infty|u_1^{\prime}|^2 r^{N+m-3}+ c_1 I_1(m,N)\int_0^\infty|u_1|^2 r^{N+m-5} dr\nonumber\\
	&\geq\sum_{k=0, k\ne 1}^{\infty} \left(\frac{N-m}{2}\right)^2\int_0^\infty|u_k^{\prime}|^2 r^{N+m-3} dr+\nonumber\\
	&+ \sum_{k=0, k\ne 1}^{\infty}c_k I_k(m,N)\int_0^\infty |u_k|^2 r^{N+m-5} dr+\nonumber\\
	&+\left(\left(\frac{N-m}{2}\right)^2-\varepsilon\right)\int_0^\infty|u_1^{\prime}|^2 r^{N+m-3}+\nonumber\\ &+\left((2c_1+\varepsilon)\left(\frac{N+m-4}{2}\right)^2\right.\nonumber\\
	&\left.+c_1^2-c_1(m-2)(m+N-4)\right)\int_0^\infty|u_1|^2 r^{N+m-5} dr.
\end{align}

The coefficient of the integral term $\int_0^\infty |u_1|^2 r^{N+m-5}dr$ in \eqref{form100}  becomes 
\begin{multline} a_1(N,m,\varepsilon):=\Bigg((2(N-1)+\varepsilon)\left(\frac{N+m-4}{2}\right)^2\\ +(N-1)^2-(N-1)(m-2)(m+N-4)\Bigg).
	\end{multline}

We choose $\varepsilon_1>0$ such that $a_1(N,m,\varepsilon_1)=0$ and we get 
\begin{equation}\label{eps}
\varepsilon_1=\frac{2(N-1)(m^2-N^2-4m+2N+2)}{(m+N-4)^2}.
\end{equation}
Indeed $\varepsilon_1>0$ because the inequality 
$ m^2-N^2-4m+2N+2 >0$ holds iff $$m\in (-\infty, 2-\sqrt{(N-1)^2+1})\cup (2+ \sqrt{(N-1)^2+1}, \infty),$$ set which contains the remaining cases $m\in (2-N, 2-\sqrt{(N-1)^2+1})$ of the case {\bf B1)}. 
From \eqref{eps} we obtain that 
$$\left(\frac{N-m}{2}\right)^2-\varepsilon_1=\frac{\left((m-2)^2-N^2\right)^2}{4(N+m-4)^2}>0.$$

	Coming back to \eqref{form100} we have
	\begin{align}
			\int_{\rr^N} |\Delta u|^2|x|^m dx &\geq \frac{\left((m-2)^2-N^2\right)^2} {4(N+m-4)^2} \sum_{k=0}^{\infty} \int_0^\infty|u_k^{\prime}|^2 r^{N+m-3} dr,
	\end{align} and therefore
	\begin{align}\label{opt4}
		\tilde{C}(N,m)\geq  \frac{\left((m-2)^2-N^2\right)^2} {4(N+m-4)^2}. 
	\end{align}

		\paragraph{Inequality \eqref{weight_ineq29} in the remaining cases of B2): $m\in ( 2+\sqrt{(N-1)^2+1}, \infty)$.}

		We want to apply the same idea as in the previous case. For any $m\in (2+\sqrt{(N+1)^2+1}$ we define the number
		\begin{equation}\label{k's}
	k(m):=\min\left\{k\in \nn \ |\ I_{k}(m,N)\geq 0 \right\}.	
		\end{equation}
 
Then
		\begin{align}
		I_1(m,N)<I_2(m,N)<...<I_{k(m)-1}<0\nonumber\\
		I_{k}(m,N)\geq I_{k(m)}(m,N)\geq 0 ,\textrm{ for all } k\geq k(m).
		\end{align}
		
		The inequality \eqref{formB2} will be rewritten as:
		\begin{align}
				\int_{\rr^N} |\Delta u|^2|x|^m dx &\geq  \left(\frac{N-m}{2}\right)^2\left(\sum_{k\geq k(m)}\int_0^\infty|u_k^{\prime}|^2 r^{N+m-3} dr+\int_0^\infty|u_0^{\prime}|^2 r^{N+m-3}dr\right) \nonumber\\
			&+ \sum_{k=1}^{k(m)-1}c_k I_k(m,N)\int_0^\infty |u_k|^2 r^{N+m-5}dr + \nonumber\\
			&=\left(\frac{N-m}{2}\right)^2\left(\sum_{k\geq k(m)}\int_0^\infty|u_k^{\prime}|^2 r^{N+m-3} dr+\int_0^\infty|u_0^{\prime}|^2 r^{N+m-3}dr\right) +\nonumber\\
			&+ \sum_{k=1}^{k(m)-1}\left(\left(\left(\frac{N-m}{2} \right)^2-\varepsilon_k\right)\int_0^\infty |u_k^{\prime}|^2 r^{N+m-3}dr \right.\nonumber\\
			&\left.+
			c_k I_k(m,N)\int_0^\infty |u_k|^2 r^{N+m-5}dr +\varepsilon_k \int_0^\infty |u_k^{\prime}|^2 r^{N+m-3}dr \right).\nonumber
			\end{align}
		Then 
		\begin{align}\label{form122}
		\int_{\rr^N} |\Delta u|^2|x|^m dx	&	\geq \left(\frac{N-m}{2}\right)^2\left(\sum_{k\geq k(m)}\int_0^\infty|u_k^{\prime}|^2 r^{N+m-3} dr+\int_0^\infty|u_0^{\prime}|^2 r^{N+m-3}dr\right) +\nonumber\\ &+\sum_{k=1}^{k(m)-1}\left(\left(\left(\frac{N-m}{2} \right)^2-\varepsilon_k\right)\int_0^\infty |u_k^{\prime}|^2 r^{N+m-3}dr +\right.\nonumber\\
			&\left.+\left(c_k I_k(m,N)+\varepsilon_k\left(\frac{N+m-4}{2}\right)^2\right)\int_0^\infty |u_k|^2 r^{N+m-5}dr \right).
		\end{align}	
We will denote the coefficient of the zero order term above by   
	\begin{equation}
		a_k(N,m,\varepsilon_k):=c_k I_k(m,N)+\varepsilon_k\left(\frac{N+m-4}{2} \right)^2.
	\end{equation}
	 We choose $\varepsilon_k$ such that $a_k(N,m,\varepsilon_k)=0$ and we obtain 
	 \begin{equation}\label{coef}
	 \varepsilon_k=\frac{-4c_k I_k(m,N)}{(N+m-4)^2}>0.
	 \end{equation}
	 For consistency, we want to make sure that  $\left(\frac{N-m}{2}\right)^2-\varepsilon_k$ is positive for every $k$, $m$ and $N$. Indeed,  
	 \begin{align}\label{coef2}
	 	\left(\frac{N-m}{2}\right)^2-\varepsilon_k&=	\left(\frac{N-m}{2}\right)^2+\frac{4c_k I_k(m,N)}{(N+m-4)^2}\\
	 	&=\frac{\left(-N+m-2k\right)^2\left(2k+m+N-4 \right)^2}{4\left(m+N-4\right)^2}>0.\nonumber
	 \end{align}
	As a consequence of \eqref{form122} we get
	 \begin{align}
	 	\int_{\rr^N} &|\Delta u|^2|x|^m dx \geq  \min_{l\leq k(m)}\left(\left(\frac{N-m}{2}\right)^2-\varepsilon_l\right)\sum_{k=0}^\infty\int_0^\infty|u_k^{\prime}|^2 r^{N+m-3}dr, 
	 \end{align}
 where $\varepsilon_l$ is as in \eqref{coef}.  
	 So, in view of \eqref{coef2} we have 
	 $$\tilde{C}(N, m)\geq \min_{l\leq k(m)} \tilde{C}(N, m, l),\quad \tilde{C}(N, m, l):=\frac{\left(-N+m-2l\right)^2\left(2l+m+N-4 \right)^2}{4\left(m+N-4\right)^2}.$$

	\subsubsection{Optimality}
In this section we aim to prove that all the lower bound constants obtained in  \eqref{opt1}, \eqref{opt2}, \eqref{opt3}, \eqref{opt4} and \eqref{coef2} are the sharp constants. For that it suffices to prove the existence of approximating sequences in \eqref{weight_ineq29} for the quoted constants. \\

\noindent {\bf Step I. The cases with radially symmetric approximations:} 
\begin{itemize}
	\item $N=1$ and $m>1$
	\item $N\geq 2$ and $m$ as in the case {\bf A)} (condition \eqref{caseA})
	\item $N\geq 2$ and $m$ as in the "good" cases {\bf B)} (condition \eqref{conditions for m})
\end{itemize}
The above cases can be treated similarly because the same sequence with radial symmetry can be built to approach   the constants \eqref{opt1}, \eqref{opt2} and \eqref{opt3}.

To prove that let us consider the radially symmetric sequence  $$u_{\epsilon}(x)=|x|^{-\frac{N+m-4}{2}+\epsilon}g(|x|)=r^{-\frac{N+m-4}{2}+\epsilon}g(r)=:U_\epsilon(r)$$
 with $g$ given in \eqref{function_g}.  
Replacing $u$ with $u_{\epsilon}$ in \eqref{form1} and arguing as in \cite{Cazacu}, since the sferical part is missing we obtain 
\begin{align}
	\int_{\mathbb{R}^N}|\Delta u_{\epsilon}|^2 |x|^m dx&=\left(\int_0^{\infty}r^{N-1+m}|U_\epsilon^{\prime\prime}(r)|^2 dr+\right.\nonumber\\
	&\left.+(N-1)(1-m)\int_0^{\infty}r^{N+m-3}|U_\epsilon^{\prime}(r)|^2dr\right)\nonumber
\end{align} 
and 
\begin{equation}
	\intr |x\cdot \n  u_{\epsilon}|^2|x|^{m-4} dx=\int_0^{\infty}|U_\epsilon^{\prime}(r)|^2 r^{N+m-3}dr.\nonumber
\end{equation}
From the definition of $U_\epsilon$ we have 
\begin{align}
	&\int_0^{\infty}  r^{N+m-3} |U^{\prime}(r)|^2dr=\nonumber\\
	& =\int_0^{\infty} r^{N+m-3}\left(\left(-\frac{N+m-4}{2}+\epsilon\right)^2   r^{-(N+m-2)+2\epsilon}g^2(r)+r^{-(N+m-4)+2\epsilon} g^{\prime}(r)^2\right)dr\nonumber\\
	&+\int_0^{\infty} r^{N+m-3}\left(2\left(-\left(\frac{N+m-4}{2}\right)+2\epsilon\right)r^{-(N+m-3)+2\epsilon}g^{\prime}(r)g(r)\right)dr\nonumber\\
	& =\frac{1}{2\epsilon}\left(-\left(\frac{N-4+m}{2}\right)+\epsilon\right)^2+\mathcal{O}(1). \label{e1}
\end{align}
Also since
\begin{align*}
	U^{\prime\prime}_{\epsilon}(r)&=\left(-\left(\frac{N-4+m}{2}\right)+\epsilon\right)\left(-\left(\frac{N-2+m}{2}\right)+\epsilon\right)r^{-\frac{N+m}{2}+\epsilon}g(r)+\mathcal{O}(1),
\end{align*}
we obtain
\begin{align}\label{e2}
	\int_0^{\infty}&r^{N+m-1}|U_\epsilon^{\prime\prime}(r)
	|^2dr=\nonumber\\
	&=\frac{1}{2\epsilon}\left(-\left(\frac{N+m-4}{2}\right)+\epsilon\right)^2
	\left(-\left(\frac{N+m-2}{2}\right)+\epsilon\right)^2+\mathcal{O}(1).
\end{align}

Due  to \eqref{e1} and \eqref{e2} we obtain
\begin{align*}
	&\frac{\intr|\Delta u_{\epsilon}|^2|x|^m dx}{\intr|x\cdot \n  u_{\epsilon}|^2|x|^{m-4} dx}=\\
	&=\frac{\frac{1}{2\epsilon}\left(\left(-\frac{N+m-4}{2}+\epsilon\right)^2\left(-\frac{N+m-2}{2}+\epsilon\right)^2+(N-1)(1-m)\left(-\frac{N+m-4}{2}+\epsilon\right)^2\right)+\mathcal{O}(1)}{\frac{1}{2\epsilon}\left(-\frac{N+m-4}{2}+\epsilon\right)^2+\mathcal{O}(1)}\\
	&=\frac{\left(-\frac{N+m-4}{2}+\epsilon\right)^2\left(\left(-\frac{N+m-2}{2}+\epsilon\right)^2+(N-1)(1-m)\right)+\mathcal{O}(\epsilon)}{\left(-\left(\frac{N+m-4}{2}\right)+\epsilon\right)^2+\mathcal{O}(\epsilon)}\\
	&\searrow\left(\frac{N+m-2}{2}\right)^2+(N-1)(1-m)=\frac{N^2-2Nm+m^2}{4}=\left(\frac{N-m}{2}\right)^2,
\end{align*} 
as $\epsilon\searrow 0$.

\noindent {\bf Step II. The "bad" cases of {\bf B)} and non radially symmetric optimal approximations:} 
\begin{itemize}
\item $N\geq 2$ and $m$ as in the "bad" cases of {\bf B1)} \\
(i.e. $m\in (2-N, 2-\sqrt{(N-1)^2+1})$)
\item $N\geq 2$ and $m$ as in the "bad" cases of {\bf B2)} \\(i.e. $m\in ( 2+\sqrt{(N-1)^2+1}, \infty)$)
\end{itemize}
In the first "bad" case  $N\geq 2$ and $m\in  (2-N,2-\sqrt{(N-1)^2+1})$
we consider the sequence  
$$	u_{\epsilon}(x)=|x|^{-\frac{N+m-4}{2}+\epsilon}g(|x|)\phi_1\left(\frac{x}{|x|}\right)=r^{-\frac{N+m-4}{2}+\epsilon}g(r)\phi_1(\sigma)=:U_\epsilon(r)\phi_1(\sigma)$$
with $g$ as in \eqref{function_g}. 

Replace once more $u$ with the above $u_{\epsilon}$ in \eqref{form1} we obtain 
\begin{align}
	\int_{\mathbb{R}^N}|\Delta u_{\epsilon}|^2 |x|^m dx&=\left(\int_0^{\infty}r^{N-1+m}|U_\epsilon^{\prime\prime}(r)|^2 dr+\right.\nonumber\\
	&+\left.\left(2c_1+(N-1)(1-m)\right)\int_0^{\infty}r^{N+m-3}|U_\epsilon^{\prime}(r)|^2dr+\right.\nonumber\\ &+\left.\left(c_1^2-c_1(m-2)(N+m-4)\right)\int_0^{\infty}r^{N+m-5}|U_\epsilon(r)|^2dr\right)\nonumber
\end{align} 
and 
\begin{equation}
	\intr |x\cdot \n  u_{\epsilon}|^2|x|^{m-4} dx=\int_0^{\infty}|U_\epsilon^{\prime}(r)|^2 r^{N+m-3}dr\nonumber.
\end{equation}

%

Due  to \eqref{e1},  \eqref{e2} and the fact that 
\begin{equation}\label{lower term}
	\int_0^{\infty}r^{N+m-5}|U_\epsilon(r)|^2dr=\frac{1}{2\epsilon}+\mathcal{O}(1),
\end{equation}
we have 
\begin{align}\label{eqst}
	\intr|\Delta u_{\epsilon}|^2|x|^m dx&=\frac{1}{2\epsilon}\left(\left(-\frac{N+m-4}{2}+\epsilon\right)^2\left(-\frac{N+m-2}{2}+\epsilon\right)^2+\right.\nonumber\\
	&+\left.\left(2c_1+(N-1)(1-m)\right)\left(-\frac{N+m-4}{2}+\epsilon\right)^2+\right.\nonumber\\
	&+\left. c_1^2-c_1(m-2)(N+m-4)\right. \Bigg)+\mathcal{O}(1)
\end{align} and 
\begin{align}\label{eqdr}
	\intr|x\cdot \n  u_{\epsilon}|^2|x|^{m-4} dx=\frac{1}{2\epsilon}\left(-\frac{N+m-4}{2}+\epsilon\right)^2+\mathcal{O}(1).
\end{align}
Due to \eqref{eqst} and \eqref{eqdr} we successively obtain
\begin{align*}
		&\frac{\intr|\Delta u_{\epsilon}|^2|x|^m dx}{\intr|x\cdot \n  u_{\epsilon}|^2|x|^{m-4} dx}=\\
		&\searrow\left(-\frac{N+m-2}{2}\right)^2+2c_1+(N-1)(1-m)+\frac{4\left(c_1^2-c_1(m-2)(N+m-4)\right)}{(N+m-4)^2}\\
		&=\frac{\left((m-2)^2-N^2\right)^2}{4\left(N+m-4\right)^2}
\end{align*}
as $\epsilon\searrow 0$.

For the second "bad" case   $N\geq 2$ and $m\in ( 2+\sqrt{(N-1)^2+1}, \infty)$
we consider the sequence  
	\begin{equation}\label{min_seq_23}
	u_\eps (x) =|x|^{-\frac{N+m-4}{2}+\eps}g(|x|)\phi_{l_{\min}}(x),
\end{equation}
where $\phi_{l_{\min}}$ is a spherical harmonic function of degree $l_{\min}$ such that\\ $\|\phi_{l_{\min}}\|_{L^2(S^{N-1})}=1$
and $l_{\min}:=\arg\min_{l\leq k(m)}\tilde{C}(N, m, l)$, where $\tilde{C}(N, m, l):=\frac{\left(-N+m-2l\right)^2\left(2l+m+N-4 \right)^2}{4\left(m+N-4\right)^2}$ and $k(m)$ was defined in \eqref{ikmn}-\eqref{k's}.

Similarly (but more technically) as above one can show that 
$$\frac{\intr|\Delta u_{\epsilon}|^2|x|^m dx}{\intr|x\cdot \n  u_{\epsilon}|^2|x|^{m-4} dx}\searrow \tilde{C}(N, m, l)$$
as $\epsilon\searrow 0$.
The details are let to the reader. 

The optimality is showed and the proof of the main result is complete now.  

\hfill

\subsection*{Acknowledgments} 
	First author was partially supported by  a grant of the Ministry of Research, Innovation and Digitization, CNCS-UEFISCDI Romania, project number PN-III-P1-1.1-TE-2021-1539, within PNCDI III.
 This work started during the elaboration of the dissertation thesis of the second author, within the master program "Advanced Studies in Mathematics" at the University of Bucharest.

\end{document}